\documentclass[12pt, reqno]{amsart}
\newcommand{\ba}{\begin{array}}
\newcommand{\ea}{\end{array}}
\newcommand{\bi}{\begin{itemize}}
\newcommand{\ei}{\end{itemize}}
\newcommand{\bc}{\begin{center}}
\newcommand{\ec}{\end{center}}
\newcommand{\bfr}{\begin{flushright}}
\newcommand{\efr}{\end{flushright}}

\theoremstyle{definition}

\numberwithin{equation}{section}

\usepackage{amssymb}
\usepackage{amsmath}

\begin{document}

\title{Convolutions with probability densities and applications to PDEs}
\author{Sorin G. Gal}

\begin{abstract}
The purpose of this paper is to introduce several new convolution operators, generated by
some known probability densities. By using the inverse Fourier transform and taking inverse steps (in the analogues of the classical procedures used for, e.g., the heat or Laplace equations), we  deduce the initial and final value problems satisfied by the new convolution integrals.
\end{abstract}
\address{University of Oradea
\newline\indent
Department of Mathematics and Computer Science
\newline\indent
Str. Universit\u{a}\c{t}ii 1
\newline\indent
410087 Oradea, Romania}
\email{galso@uoradea.ro}

\maketitle



\section{Introduction}

It is well known the fact that the classical Gauss-Weierstrass, Poisson-Cauchy and Picard convolution singular integrals are based on convolutions with the standard normal density function $\frac{1}{\sqrt{\pi}}\cdot e^{-x^{2}}$, standard Cauchy density function $\frac{1}{\pi}\cdot \frac{1}{1+x^{2}}$ and Laplace density function $\frac{1}{2}\cdot e^{-|x|}$, respectively. Their approximation properties are studied, for example, in \cite{Butzer}, \cite{Gal1}. Also, by using the Fourier transform method, it is known that the solutions of the initial value problems for the heat equation and Laplace equation are exactly the Gauss-Weierstrass and Poisson-Cauchy convolution singular integrals, respectively, see, e.g., \cite{Gold}, p. 23. On the other hand, in our best knowledge, the initial value problem and the partial differential equation corresponding to the Picard singular integral, is missing from mathematical literature.
The main aim of the present paper is somehow inverse : introducing convolution singular integrals based on some known probability densities, we use the inverse Fourier transform in order to find the partial differential equations (initial and final value problems) satisfied by these integrals, including the Picard singular integral.

\section{Definitions of Convolution Operators}

In this section we introduce several convolution operators, based on some well-known densities of probability.

If $d(t, x)\ge 0$ with $t>0$ and $x\in \mathbb{R}$ is a probability density, that is $\int_{-\infty}^{+\infty}d(t, x)d x=1$, then our definitions are based on the general known formula
$$O_{t}(f)(x)$$
\begin{equation}\label{eq0}
=d(t, \cdot)*f(\cdot)=\int_{-\infty}^{+\infty}f(u)\cdot d(t, x-u) du=\int_{-\infty}^{+\infty}f(x-v)\cdot d(t, v) dv.
\end{equation}

{\bf Definition 2.1} (i) For the Maxwell-Boltzmann type probability density (see, e.g., \cite{Papou}, p. 104 and pp. 148-149)
$$d(t, x)=\frac{1}{\sqrt{2 \pi}}\cdot  \frac{x^{2} e^{-x^{2}/(2 t^{2})}}{t^{3}}, x\in \mathbb{R}, t>0$$
and $f:\mathbb{R}\to \mathbb{R}$, we can formally define the Maxwell-Boltzmann convolution operator
\begin{equation}\label{eq17}
S_{t}(f)(x)=\frac{1}{\sqrt{2 \pi}}\int_{-\infty}^{+\infty}f(x-v)\cdot \frac{v^{2} e^{-v^{2}/(2 t^{2})}}{t^{3}} d v, t>0, x\in \mathbb{R}.
\end{equation}

(ii) For the Laplace type probability density (see, e.g., \cite{Ever})
$$d(t, x)=\frac{1}{2 t}e^{-|x|/t}, t>0, x\in \mathbb{R}$$
and $f:\mathbb{R}\to \mathbb{R}$, we can formally define the classical Picard convolution operator
\begin{equation}\label{eq24}
P_{t}(f)(x)=\frac{1}{2 t}\int_{-\infty}^{+\infty}f(x-v)\cdot e^{-|v|/t} d v, t>0, x\in \mathbb{R}.
\end{equation}

(iii) For the exponential probability density (see, e.g., \cite{Ever}, \cite{John})
$$d(t, x)=\frac{t e^{-t|x|}}{2}, x\in \mathbb{R}, t>0$$
and $f:\mathbb{R}\to \mathbb{R}$, we can formally define the exponential convolution operator
\begin{equation}\label{eq8}
E_{t}(f)(x)=\int_{-\infty}^{+\infty}f(x-v)\cdot \frac{t e^{-t|v|}}{2} d v, t>0, x\in \mathbb{R}.
\end{equation}

(iv) For any $n\in \mathbb{N}$, $P_{t}(f)(x)$ can be generalized to the so called Jackson type generalization of the Picard singular integral defined by (see, e.g., \cite{Gal1})
$$P_{n, t}(f)(x)=-\frac{1}{2 t}\int_{-\infty}^{+\infty}\left (\sum_{k=1}^{n+1}(-1)^{k}{n+1 \choose k}f(x+k v)e^{-|v|/t}\right )d v$$
$$=\int_{-\infty}^{+\infty}f(x-u)\left [\sum_{k=1}^{n+1}(-1)^{k+1}{n+1 \choose k}\cdot \frac{1}{k}\cdot \frac{e^{-|u|/(kt)}}{2 t}\right ]du, t>0, x\in \mathbb{R}.$$

(v) Starting from the well known Gauss-Weierstrass operator $W_{t}(f)(x)=\frac{1}{\sqrt{ \pi t}}\int_{-\infty}^{+\infty}f(x-v)e^{-v^2/t}d v$, we can define its Jackson type generalization by (see, e.g., \cite{Gal1})
$$W_{n, t}(f)(x)=-\frac{1}{2 C^{*}(t)}\cdot \int_{-\infty}^{+\infty}\left (\sum_{k=1}^{n+1}(-1)^{k}{n+1 \choose k}f(x+k v)e^{-v^2/t}\right )d v$$
$$=\int_{-\infty}^{+\infty}f(x-u)\left [\sum_{k=1}^{n+1}(-1)^{k+1}{n+1 \choose k}\cdot \frac{1}{k}\cdot \frac{e^{-u^{2}/(k t)}}{2 C^{*}(t)}\right ]du, t>0, x\in \mathbb{R},$$
where $C^{*}(t)=\int_{0}^{\infty}e^{-u^{2}/t}du=\frac{\sqrt{t \pi}}{2}.$
Therefore, for any $n\in \mathbb{N}$, we can write
$$W_{n, t}(f)(x)=\int_{-\infty}^{+\infty}f(x-u)\left [\sum_{k=1}^{n+1}(-1)^{k+1}{n+1 \choose k}\cdot \frac{1}{k}\cdot \frac{e^{-u^{2}/(k t)}}{\sqrt{\pi t}}\right ]du.$$

\section{Applications to PDE}

Concerning the convolution operators defined in Section 2, we can state the following applications to PDE.

{\bf Theorem 3.1.} {\it (i) Suppose that $f, f^{\prime}, f^{\prime \prime}, f^{\prime \prime \prime}, f^{(4)}:\mathbb{R}\to \mathbb{R}$ are bounded and uniformly continuous on $\mathbb{R}$.
The solution of the initial value problem
$$\frac{\partial u}{\partial t}(x, t)=
t^3 \frac{\partial^4 u}{\partial x^4}(x, t)-t^2 \frac{\partial^3 u}{\partial x^2 \partial t}(x, t)+3t \frac{\partial^2 u}{\partial x^2}(x, t),$$
$$\lim_{s\searrow 0}u(x, s)=f(x), t>0, x\in \mathbb{R},$$
is $u(x, t):=S_{t}(f)(x)$.

(ii) Suppose that $f, f^{\prime}, f^{\prime \prime} :\mathbb{R}\to \mathbb{R}$ are bounded and uniformly continuous on $\mathbb{R}$. The solution of the initial value problem
$$\frac{\partial u}{\partial t}(x, t)=t^{2} \frac{\partial^3 u}{\partial x^2 \partial t}(x, t)+2t \frac{\partial^2 u}{\partial x^2}(x, t),\, \, \, \lim_{s\searrow 0}u(x, s)=f(x), t>0, x\in \mathbb{R}$$
is $u(x, t):=P_{t}(f)(x)$.

(iii) Suppose that $f, f^{\prime}, f^{\prime \prime} :\mathbb{R}\to \mathbb{R}$ are bounded and uniformly continuous on $\mathbb{R}$.
The solution of the final value problem
$$\frac{\partial u}{\partial t}(x, t)=\frac{1}{t^{2}}\cdot \frac{\partial^{3} u}{\partial x^2 \partial t}(x, t)-
\frac{2}{t^{3}}\cdot \frac{\partial^{2} u}{\partial x^{2}}(x, t), \, \, \, \lim_{s\to \infty}u(x, s)=f(x), t>0, x\in \mathbb{R}$$
is $u(x, t):=E_{t}(f)(x)$.

(iv) Suppose that $f, f^{\prime}, f^{\prime \prime} :\mathbb{R}\to \mathbb{R}$ are bounded and uniformly continuous on $\mathbb{R}$.
We have
$$P_{n, t}(f)(x)=\sum_{k=1}^{n+1}(-1)^{k+1}{n+1 \choose k}\cdot u_{k}(x, t),$$
where $u_{k}(x, t)=P_{k t}(f)(x)=\frac{1}{2 k t}\cdot \int_{-\infty}^{+\infty}f(x-u)e^{-|u|/(k t)}d u$, $k=1, ..., n+1$ are solutions of the initial value problems (for $t>0$ and $x\in \mathbb{R}$)
$$\frac{\partial u_{k}}{\partial t}(x, t)=k^{2}t^{2}\frac{\partial^{3} u_{k}}{\partial x^{2}\partial t}(x, t)+2k^{2}t\cdot \frac{\partial^{2} u_{k}}{\partial x^{2}}(x, t), \, \, \, \lim_{s\searrow 0}u_{k}(x, s)=f(x).$$

(v) Suppose that $f, f^{\prime}, f^{\prime \prime} :\mathbb{R}\to \mathbb{R}$ are bounded and uniformly continuous on $\mathbb{R}$. We have
$$W_{n, t}(f)(x)=\sum_{k=1}^{n+1}(-1)^{k+1}{n+1 \choose k}\cdot u_{k}(x, t),$$
where $u_{k}(x, t)=\frac{1}{\sqrt{k}}W_{k t}(f)(x)=\frac{1}{k\sqrt{\pi t}}\cdot \int_{-\infty}^{+\infty}f(x-u)e^{-u^{2}/(k t)}d u$, $k=1, ..., n+1$ are solutions of the initial value problems
$$\frac{\partial u_{k}}{\partial t}(x, t)=\frac{k}{4}\cdot \frac{\partial^{2} u_{k}}{\partial x^{2}}(x, t), \, \, \, \lim_{s\searrow 0}u_{k}(x, s)=\frac{1}{\sqrt{k}}f(x), t>0, x\in \mathbb{R}.$$}

{\bf Proof.} Since for the convolution operator given by (\ref{eq0}), in general we have $d(t, v)\ge 0$, for all $t>0$ and $v\in \mathbb{R}$, by the standard method we easily get
$$|O_{t}(f)(x)-f(x)|\le \int_{-\infty}^{+\infty}|f(x-v)-f(x)|d(t, v) d v$$
$$\le \int_{-\infty}^{+\infty}\omega_{1}(f; |v|)_{\mathbb{R}}d(t, v) d v
\le 2\omega_{1}(f ; \varphi(t))_{\mathbb{R}},$$
where $\omega_{1}(f ; \delta)_{\mathbb{R}}=\sup\{|f(x)-f(y)|; x, y \in \mathbb{R}, |x-y|\le \delta\}$ and $\varphi(t)=\int_{-\infty}^{+\infty}|v|\cdot d(t, v)d v$.

Evidently that this method is useful only if $\varphi(t)<+\infty$ for all $t>0$.

In order to deduce the PDE equations satisfied by various convolution operators, we will need the concepts of Fourier transform of a function $g$, defined by
$$F(g)(\xi)=\hat{g}(\xi)=\frac{1}{\sqrt{2\pi}}\cdot \int_{-\infty}^{+\infty}g(x)e^{-i \xi x}d x, \mbox{ if } \int_{-\infty}^{+\infty}|g(x)|d x<+\infty,$$
and of inverse Fourier transform defined by
$$F^{-1}(\hat{g})(x)=g(x)=\frac{1}{\sqrt{2 \pi}}\cdot\int_{-\infty}^{+\infty}\hat{g}(\xi)e^{i \xi x}d \xi.$$

(i) By making the change of variable $v=\sqrt{2} t s$, we get
$$\varphi(t)=\frac{1}{t^{3}}\cdot \frac{\sqrt{2}}{\sqrt{\pi}}\int_{0}^{\infty}v^{3}e^{-v^{2}/(2 t^{2})}dv
=\frac{1}{t^{3}}\cdot \frac{\sqrt{2}}{\sqrt{\pi}}\int_{0}^{\infty}(2\sqrt{2}t^{3} s^{3})e^{-s^{2}}(\sqrt{2} t )d s$$
$$=\frac{4\sqrt{2}}{\sqrt{\pi}} t\int_{0}^{\infty}s^{3}e^{-s^{2}}d s=\frac{2\sqrt{2}}{\sqrt{\pi}} t <2 t,$$
which immediately implies
$$|S_{t}(f)(x)-f(x)|\le 4\omega_{1}(f; t)_{\mathbb{R}}, t>0, x\in \mathbb{R}.$$
Taking into account the uniform continuity of $f$, the above inequality immediately implies that $\lim_{t\searrow 0}S_{t}(f)(x)=f(x)$, for all $x\in \mathbb{R}$. Therefore we may take, by convention, $S_{0}(f)(x)=f(x)$, for all $x\in \mathbb{R}$.

Now, in order to deduce the PDE satisfied by $S_{t}(f)(x)$, we write it in the form
$$S_{t}(f)(x)$$
$$=\frac{1}{\sqrt{2 \pi}}\int_{-\infty}^{+\infty}f(y)\cdot \frac{(x-y)^{2} e^{-(x-y)^{2}/(2 t^{2})}}{t^{3}} d y=\frac{1}{\sqrt{2 \pi}}\int_{-\infty}^{+\infty}f(y)\cdot \hat{g}_{t}(y-x)d y.$$
Here, by using standard reasonings/calculation (or the WolframAlpha soft of calculation), we obtain
$$g_{t}(\xi)=F^{-1}_{w}[w^{2}e^{-w^{2}/(2 t^{2})}/t^{3}](\xi, t)=e^{-t^{2}\xi^{2}/2}(1-t^{2}\xi^{2}),$$
which implies
$$S_{t}(f)(x)=\frac{1}{\sqrt{2 \pi}}\cdot \int_{-\infty}^{+\infty}f(y)\left [\frac{1}{\sqrt{2 \pi}}\int_{-\infty}^{+\infty}e^{-i(y-x)\xi}e^{-t^{2}\xi^2/2}(1-t^{2}\xi^{2})d \xi\right ]d y$$
$$=\frac{1}{\sqrt{2 \pi}}\int_{-\infty}^{+\infty}\left [\frac{1}{\sqrt{2 \pi}}\cdot \int_{-\infty}^{+\infty} e^{-i y \xi}f(y)dy\right ]e^{i x \xi} e^{-t^{2}\xi^2/2}(1-t^{2}\xi^{2})d \xi$$
$$=\frac{1}{\sqrt{2 \pi}}\cdot \int_{-\infty}^{+\infty}e^{i x \xi}\hat{f}(\xi)e^{-t^{2}\xi^2/2}(1-t^{2}\xi^{2})d \xi$$
$$=\frac{1}{\sqrt{2 \pi}}\cdot \int_{-\infty}^{+\infty}e^{i x \xi}\hat{u}(\xi, t)d \xi:=u(x, t),$$
where
$$\hat{u}(\xi, t)=\hat{f}(\xi)\cdot e^{-t^{2}\xi^2/2}(1-t^{2}\xi^{2}).$$
This is equivalent to $\hat{u}(\xi, t)\cdot \frac{e^{t^{2}\xi^2/2}}{1-t^{2}\xi^{2}}=\hat{f}(\xi)$, which is equivalent to
$$\frac{\partial }{\partial t}\left [\hat{u}(\xi, t)\cdot \frac{e^{t^{2}\xi^2/2}}{1-t^{2}\xi^{2}}\right ]=
\frac{\partial \hat{u}}{\partial t}(\xi, t)\cdot \frac{e^{t^2 \xi^2/2}}{1-t^2 \xi^2}+\hat{u}(\xi, t)\cdot \left (\frac{e^{t^2 \xi^2/2}}{1-t^2 \xi^2}\right )^{\prime}_{t}=0.$$
Note that the above relation can be evidently written under the form
$$\frac{\partial }{\partial t}\left [\hat{u}(\xi, t)\cdot \frac{1}{F^{-1}_{w}(d(t, w))(\xi, t)}\right ]=0,$$
where $d(t, x)$ is the Maxwell-Boltzmann type probability density in Definition 2.1, (i), entering in the formula for $S_{t}(f)(x)$.

After simple calculation, the above formula is formally equivalent to (of course for $1\not=t^2 \xi^2$)
$$\frac{\partial \hat{u}}{\partial t}(\xi, t)+t^{2}\left (-\xi^{2}\cdot \frac{\partial \hat{u}}{\partial t}(\xi, t)\right )-3t[-\xi^{2}\hat{u}(\xi, t)]-t^3\cdot [\xi^4 \hat{u}(\xi, t)]=0.$$
Now, taking into account that
$$\frac{\partial \hat{u}}{\partial t}(\xi, t)=\widehat{\frac{\partial u}{\partial t}}(\xi, t), \,
\widehat{\frac{\partial^{2} u}{\partial x^2}}(\xi, t) = -\xi^{2}\hat{u}(\xi, t), \, \widehat{\frac{\partial^{4} u}{\partial x^4}}(\xi, t) = \xi^{4}\hat{u}(\xi, t),$$
and replacing above, we obtain
$$F\left (\frac{\partial u}{\partial t}+t^2 \frac{\partial^3 u}{\partial x^2 \partial t}-3t \frac{\partial^2 u}{\partial x^2}-t^3 \frac{\partial^4 u}{\partial x^4}\right )(\xi, t)=0,$$
that is
$$\frac{\partial u}{\partial t}(x, t)=
t^3 \frac{\partial^4 u}{\partial x^4}(x, t)-t^2 \frac{\partial^3 u}{\partial x^2 \partial t}(x, t)+3t \frac{\partial^2 u}{\partial x^2}(x, t).$$
Finally, following the above steps in inverse order, we arrive at the conclusion in the statement.

(ii) By \cite{Butzer}, p. 142, Corollary 3.4.2, it was obtained
$$|f(x)-P_{t}(f)(x)|\le C\omega_{2}(f; t)_{\mathbb{R}}.$$
Therefore, it is immediate that $\lim_{t\searrow 0}P_{t}(f)(x)=f(x)$, for all $x\in \mathbb{R}$.

In order to deduce the PDE satisfied by $P_{t}(f)(x)$, we reason exactly as in the above case (i). Indeed, by standard calculation (or by making use of the WolframAlpha program), we get
$$F^{-1}_{w}[e^{-|w|/t}/(2t)](\xi, t)=\frac{1}{\sqrt{2 \pi}}\cdot \frac{1}{1+t^2 \xi^2}$$
and similar reasonings with those in the case (i), immediately leads to
$$P_{t}(f)(x)=\frac{1}{\sqrt{2 \pi}}\int_{-\infty}^{+\infty}e^{i x \xi}\hat{u}(\xi, t)d \xi:=u(x, t),$$
where
$$\hat{u}(\xi, t)=\hat{f}(\xi)\cdot \frac{1}{t^2 \xi^2 + 1}.$$
In fact, directly as in the proof of Theorem 3.1, (i), we can write
$$\frac{\partial }{\partial t}\left [\hat{u}(\xi, t)\cdot \frac{1}{F^{-1}_{w}(d(t, w))(\xi, t)}\right ]=0,$$
where $d(t, x)$ is the Laplace type probability density in Definition 2.1, (ii), entering in the formula for $P_{t}(f)(x)$.

Therefore,
$$\frac{\partial }{\partial t}\left [\hat{u}(\xi, t)\cdot (1+t^{2}\xi^{2})\right ]=
\frac{\partial \hat{u}}{\partial t}(\xi, t)+t^2 \xi^{2}\frac{\partial \hat{u}}{\partial t}(\xi, t)+2t \xi^2 \hat{u}(\xi, t)=0,$$
which immediately leads to
$$\frac{\partial u}{\partial t}(x, t)=t^{2} \frac{\partial^3 u}{\partial x^2 \partial t}(x, t)+2t \frac{\partial^2 u}{\partial x^2}(x, t).$$
Following the above steps, now from the end to the beginning, we arrive at the conclusion in the statement.

(iii) Firstly, we observe that $E_{t}(f)(x)=P_{1/t}(f)(x)$, for all $t>0$ and $x\in \mathbb{R}$. The, by (ii) we immediately get
$$|E_{t}(f)(x)-f(x)|=|P_{1/t}(f)(x)-f(x)|\le 2\omega_{1}\left (f; \frac{1}{t}\right )_{\mathbb{R}}.$$
Then, again by standard calculation (or by using WolframAlpha), we have  $F^{-1}_{w}(e^{-|w|t})(\xi, t)=\frac{\sqrt{2}}{\sqrt{\pi}}\cdot \frac{t}{t^{2}+\xi^{2}}$. This immediately implies
$$F^{-1}(d(t, w))(\xi, t)=\frac{t}{2}\cdot F^{-1}_{w}(e^{-|w|t})(\xi, t)=\frac{1}{\sqrt{2 \pi}}\cdot \frac{t^2}{t^{2}+\xi^{2}}.$$
It follows $\frac{1}{F^{-1}(d(t, w))(\xi, t)}=\sqrt{2 \pi}\cdot \frac{t^{2}+\xi^{2}}{t^{2}}=\sqrt{2 \pi}\left (1+\frac{\xi^{2}}{t^{2}}\right )$.

Therefore, denoting $u(x, t)=E_{t}(f)(x)$, by the  method used at the above points, we arrive at the PDE
$$\frac{\partial}{\partial t}\left (\hat{u}(\xi, t)\cdot \left (1+\frac{\xi^{2}}{t^{2}} \right )\right )=
\frac{\hat{u}}{\partial t}(\xi, t)\left (1+\frac{\xi^{2}}{t^{2}}\right )+\hat{u}(\xi, t)\left (-\frac{2\xi^{2}}{t^{3}}\right )=0.$$
This immediately leads to the following PDE, satisfied by $u(x, t)=E_{t}(f)(x)$
$$\frac{\partial u}{\partial t}(x, t)=\frac{1}{t^{2}}\cdot \frac{\partial^{3} u}{\partial x^2 \partial t}(x, t)-
\frac{2}{t^{3}}\cdot \frac{\partial^{2} u}{\partial x^{2}}(x, t).$$

Since $E_{t}(f)(x)=P_{1/t}(f)(x)$, it follows that $\lim_{t\nearrow \infty}E_{t}(f)(x)=f(x)$, for all $x\in \mathbb{R}$.

Following the above steps in inverse order, we arrive at the conclusion in the statement.

(iv) Concerning the approximation properties of $P_{n, t}(f)(x)$, in \cite{Gal1} it was obtained the estimate
$$|f(x)-P_{n, t}(f)(x)|\le \sum_{k=1}^{n+1} k! \cdot {n+1 \choose k}\cdot \omega_{n+1}(f ; t)_{\mathbb{R}},$$
where $\omega_{n+1}(f ; \delta)=\sup_{0\le h\le \delta}\{|\Delta_{h}^{n+1}f(x) ; x\in \mathbb{R}\}$, with
$\Delta_{h}^{n+1}=\sum_{j=0}^{n+1}(-1)^{n+1-j}{n+1 \choose j}f(x+j h)$. This immediately implies $\lim_{t\searrow 0}P_{n, t}(f)(x)=f(x)$, for all $x\in \mathbb{R}$.

In order to deduce the PDE satisfied by $P_{n, t}(f)(x)$, since $F^{-1}_{w}$ is linear operator and since known calculation (or by using the WolframAlpha software) give
$F^{-1}_{w}(e^{-|w|/t})(\xi, t)=\frac{\sqrt{2}}{\sqrt{\pi}}\cdot \frac{t}{t^{2}\xi^{2}+1}$, replacing here $t$ by $k t$, we easily obtain
$$F^{-1}_{w}\left [\sum_{k=1}^{n+1}(-1)^{k+1}{n+1 \choose k}\cdot \frac{1}{2 t k}e^{-|w|/(k t)}\right ](\xi, t)$$
$$=\frac{1}{\sqrt{2 \pi}}\cdot \sum_{k=1}^{n+1}(-1)^{k+1}{n+1 \choose k}\cdot \frac{1}{k t}\cdot \frac{k t}{k^{2}t^{2}\xi^{2}+1}$$
$$=\frac{1}{\sqrt{2 \pi}}\cdot \sum_{k=1}^{n+1}(-1)^{k+1}{n+1 \choose k}\cdot \frac{1}{k^{2}t^{2}\xi^{2}+1}.$$
Therefore, denoting $u(x, t):=P_{n, t}(f)(x)$ we immediately get the differential equation
$$\frac{\partial }{\partial t}\left [\hat{u}(\xi, t)\cdot \frac{1}{\sum_{k=1}^{n+1}(-1)^{k+1}{n+1 \choose k}\cdot \frac{1}{k^{2}t^{2}\xi^{2}+1}}\right ]=0,$$
which is equivalent to
$$\sum_{k=1}^{n+1}(-1)^{k+1}{n+1 \choose k}\left [\frac{\partial \hat{u}}{\partial t}(\xi, t)\cdot \frac{1}{k^{2}t^{2}\xi^{2}+1}+\hat{u}(\xi, t)\cdot \frac{2 k^{2} t \xi^{2}}{(k^{2}t^{2}\xi^{2}+1)^{2}}\right ]=0.$$
It is worth noting that denoting
$$u_{k}(x, t)=P_{k t}(f)(x)=\frac{1}{2 k t}\cdot \int_{-\infty}^{+\infty}f(x-u)e^{-|u|/(k t)}d u,$$
we can write
$$P_{n, t}(f)(x)=\sum_{k=1}^{n+1}(-1)^{k+1}{n+1 \choose k}\cdot u_{k}(x, t),$$
where reasoning as above for $P_{t}(f)(x)$, we easily obtain
$$\frac{\partial \widehat{u_{k}}}{\partial t}(\xi, t)=k^{2}t^{2}\frac{\partial^{3} \widehat{u_{k}}}{\partial x^{2}\partial t}(\xi, t)+2k^{2}t\cdot \frac{\partial^{2}\widehat{u_{k}}}{\partial x^{2}}(\xi, t)$$
and which implies
$$\frac{\partial u_{k}}{\partial t}(x, t)=k^{2}t^{2}\frac{\partial^{3} u_{k}}{\partial x^{2}\partial t}(x, t)+2k^{2}t\cdot \frac{\partial^{2} u_{k}}{\partial x^{2}}(x, t), \, x\in \mathbb{R}, t>0, k=1, ..., n+1,$$
with $u_{k}(x, 0)=f(x)$, for all $x\in \mathbb{R}$, $k=1, ..., n+1$.

(v) Concerning the approximation properties of $W_{n, t}(f)(x)$, reasoning as in \cite{Gal1}, we get the estimate
$$|f(x)-W_{n, t}(f)(x)|\le C_{n}\cdot \omega_{n+1}(f ; \sqrt{t})_{\mathbb{R}},$$
where $C_{n}>0$ is a constant independent of $f$, $t$ and $x$. This immediately implies that $\lim_{t\searrow 0}W_{n, t}(f)(x)=f(x)$, for all $x\in \mathbb{R}$.

Now, in order to deduce the PDE satisfied by $W_{n, t}(f)(x)$, since $F^{-1}_{w}$ is linear operator and since (by, e.g., WolframAlpha software) we have $F^{-1}_{w}(e^{-w^2/t})(\xi, t)=\frac{\sqrt{t}\cdot e^{-t \xi^{2}/4}}{\sqrt{2}}$, replacing here $t$ by $k t$, we easily obtain $F^{-1}_{w}(e^{-w^{2}/(k t)})(\xi, t)=\frac{\sqrt{k t} e^{-k t\xi^{2}/4}}{\sqrt{2}}$ and
$$F^{-1}_{w}\left [\sum_{k=1}^{n+1}(-1)^{k+1}{n+1 \choose k}\cdot \frac{1}{\sqrt{\pi t} k}e^{-w^{2}/(k t)}\right ](\xi, t)$$
$$=\frac{1}{\sqrt{\pi}}\cdot \sum_{k=1}^{n+1}(-1)^{k+1}{n+1 \choose k}\cdot \frac{1}{k \sqrt{t}}\cdot \frac{\sqrt{k t} e^{-kt\xi^{2}/4}}{\sqrt{2}}$$
$$=\frac{1}{\sqrt{2 \pi}}\cdot \sum_{k=1}^{n+1}(-1)^{k+1}{n+1 \choose k}\cdot \frac{1}{\sqrt{k}}\cdot e^{-k t\xi^{2}/4}.$$
Therefore, denoting $u(x, t):=W_{n, t}(f)(x)$ we immediately get the differential equation
$$\frac{\partial }{\partial t}\left [\hat{u}(\xi, t)\cdot \frac{1}{\sum_{k=1}^{n+1}(-1)^{k+1}{n+1 \choose k}\cdot \frac{1}{\sqrt{k}}\cdot e^{-kt\xi^{2}/4}}\right ]=0,$$
which is equivalent to
$$\sum_{k=1}^{n+1}(-1)^{k+1}{n+1 \choose k}\cdot \frac{1}{\sqrt{k}}\left [\frac{\partial \hat{u}}{\partial t}(\xi, t)\cdot e^{-kt\xi^{2}/4}+\hat{u}(\xi, t)\cdot \frac{k \xi^{2}}{4}\cdot e^{-kt\xi^{2}/4}\right ]=0.$$
It is worth noting that denoting
$$u_{k}(x, t)=\frac{1}{\sqrt{k}}W_{k t}(f)(x)=\frac{1}{k \sqrt{\pi t}}\cdot \int_{-\infty}^{+\infty}f(x-u)e^{-u^{2}/(k t)}d u,$$
we can write
$$W_{n, t}(f)(x)=\sum_{k=1}^{n+1}(-1)^{k+1}{n+1 \choose k}\cdot u_{k}(x, t),$$
where reasoning as above but now for $W_{t}(f)(x)$, we easily obtain
$$\frac{\partial \widehat{u_{k}}}{\partial t}(\xi, t)=\frac{k}{4}\cdot \frac{\partial^{2} \widehat{u_{k}}}{\partial x^{2}}(\xi, t)$$
and which implies
$$\frac{\partial u_{k}}{\partial t}(x, t)=\frac{k}{4}\cdot \frac{\partial^{2} u_{k}}{\partial x^{2}}(x, t), \, x\in \mathbb{R}, t>0, k=1, ..., n+1,$$
with $u_{k}(x, 0)=f(x)$, for all $x\in \mathbb{R}$, $k=1, ..., n+1$. $\hfill \square$

{\bf Remark 3.2} The methods in this paper could be used to make analogous studies for the convolutions with other
known probability densities, like the Rayleigh probability density (shortly written Rayleigh p.d.), Gumbel p.d., logistic p.d., Johnson p.d., Fr\'echet p.d., Gompetz p.d., L\'evy p.d., Lomax p.d. and so on.

{\bf Remark 3.3.} It would be also of interest to use the methods in this paper to the case of the corresponding complex convolutions, based on the ideas and results in the books \cite{Gal2} and \cite{GGG}.

\end{document}